\documentclass[11pt,letterpaper]{article}
\usepackage[T1]{fontenc}
\usepackage{lmodern}
\usepackage[margin=1in]{geometry}
\usepackage{amsmath,amssymb,amsthm}
\usepackage{microtype}
\usepackage[hidelinks]{hyperref}

\newtheorem{theorem}{Theorem}[section]
\newtheorem{proposition}[theorem]{Proposition}
\newtheorem{lemma}[theorem]{Lemma}
\newtheorem{corollary}[theorem]{Corollary}
\theoremstyle{remark}
\newtheorem{remark}[theorem]{Remark}

\newcommand{\shorttitle}{Sharp fractional clique packings}
\title{Sharp fractional clique packings in clique-free graphs}
\author{Jialin He\thanks{School of Mathematical Sciences, Key Laboratory of MEA, and Shanghai Key Laboratory of PMMP, East China Normal University, Shanghai 200241, China. Emails: \href{mailto:jlhe@math.ecnu.edu.cn}{jlhe@math.ecnu.edu.cn}. Supported in part by the National Natural Science Foundation of China (Grant No.\ 12601663) and the Science and Technology Commission of Shanghai Municipality (Grant No.\ 22DZ2229014).}
  \and Fuhong Ma\thanks{Corresponding author. School of Mathematics and Statistics, Shandong University of Technology, Zibo 255000, China. Email: \href{mailto:mafuhongsdnu@163.com}{mafuhongsdnu@163.com}. Supported by the National Natural Science Foundation of China (Grant No.\ 12401458) and the Natural Science Foundation of Shandong Province (Grant No.\ ZR2024QA239).}}
\date{}

\begin{document}
\maketitle
\markboth{Jialin He and Fuhong Ma}{\shorttitle}

\begin{abstract}
We study fractional clique packings in graphs with forbidden cliques,
focusing on the optimal lower bound in terms of the excess
$e(G)-\rho_{r-1}n^2$ above the asymptotic Tur\'an benchmark.
For integers $r\ge3$ and $M\ge r$, we prove that every $n$-vertex
$K_{M+1}$-free graph $G$ satisfies
\[
\nu_r^*(G)\ge
\gamma_{r,M}\bigl(e(G)-\rho_{r-1}n^2\bigr),
\qquad
\gamma_{r,M}:=
\min\left\{1,\frac{2(M-1)}{r(M-r+1)}\right\},
\]
and show that $\gamma_{r,M}$ is best possible for every pair $(r,M)$.
Here, $\nu_r^*(G)$ denotes the maximum total weight of a fractional
$K_r$-packing in $G$, and
$\rho_{r-1}=(r-2)/(2(r-1))$.
We further determine the optimal asymptotic coefficient for integer
$K_r$-packings in $H$-free graphs, for every fixed graph $H$ with
$\chi(H)\ge r+1$, uniformly over all positive excesses.
In particular, this coefficient equals $1$ if and only if
$\chi(H)\le r+2$. Our proof combines fractional symmetrization,
positive semidefinite matrix methods, and a sampling argument derived
from a sharp approximate Carath\'eodory theorem. These results
quantify how a forbidden-clique condition strengthens fractional
clique-packing bounds.
\end{abstract}


\section{Introduction}\label{sec:introduction}

Tur\'an's theorem gives the edge threshold for forcing a clique.
The corresponding packing problem asks how many edge-disjoint cliques
are forced by the excess above this threshold.
For arbitrary graphs, the best uniform asymptotic coefficient is known.
We determine how this coefficient improves when a clique of a prescribed
larger order is forbidden.

All graphs are finite and simple. For $r\ge3$, let $\nu_r(G)$ be the
maximum number of pairwise edge-disjoint copies of $K_r$ in $G$.
A fractional $K_r$-packing assigns nonnegative weights to copies of
$K_r$ so that the total weight on each edge is at most one; its maximum
total weight is denoted by $\nu_r^*(G)$. Write $T_j(n)$ for the balanced
complete $j$-partite graph on $n$ vertices and $t_j(n)=e(T_j(n))$.
We put $\rho_j=(j-1)/(2j)$, so $t_j(n)=\rho_j n^2-O_j(1)$.
An $H$-free graph contains no subgraph isomorphic to $H$; the
subgraph is not required to be induced.

The packing problem is closely related to classical decomposition questions.
Erd\H{o}s, Goodman, and P\'osa~\cite{EGP} proved that every $n$-vertex
graph can be decomposed into at most $t_2(n)$ edges and triangles.
If such a decomposition contains $a$ triangles, then it has
$e(G)-2a$ members. Consequently,
$\nu_3(G)\ge(e(G)-t_2(n))/2$.
Kr\'al', Lidick\'y, Martins, and Pehova~\cite{KLMP} obtained an
asymptotic lower bound of $n^2/2$ for the minimum total order of the
members of an edge-and-triangle decomposition. Blumenthal, Lidick\'y,
Pehova, Pfender, Pikhurko, and Volec~\cite{BLPPPV} later determined
the sharp bound and its extremal graphs for sufficiently large $n$.
Their result implies
$\nu_3(G)\ge(2/3)(e(G)-t_2(n))-O(1)$.

For fixed $r\ge3$, Gy\H{o}ri~\cite[Conjecture~2]{Gyori1992}
conjectured that every graph with $t_{r-1}(n)+k$ edges contains at
least $2k/r-o(n^2)$ edge-disjoint copies of $K_r$.
Balogh and Wigal~\cite{BW} proved the conjecture in the relative-error
form $\nu_r(G)\ge(2/r-o(1))k$ for $k>0$.
They established the finite fractional inequality
\begin{equation}\label{eq:bw}
\nu_r^*(G)\ge\frac2r\bigl(e(G)-\rho_{r-1}n^2\bigr).
\end{equation}
Their symmetrization reduces the fractional inequality to complete
multipartite graphs. To obtain the integer-packing bound, they combine
Haxell and R\"odl's theorem~\cite{HR}, in the uniform form stated as
Theorem~1.7 of Balogh and Wigal~\cite{BW}, which gives
$\nu_r^*(G)-\nu_r(G)=o(n^2)$ uniformly over $n$-vertex graphs, with
Gy\H{o}ri's small-excess estimate~\cite{Gyori1991,Gyori1992}.
The latter is stated for every fixed $r\ge3$ in Theorem~1.6 of
Balogh and Wigal~\cite{BW}:
$\nu_r(G)\ge k-O_r(k^2/n^2)$ when $k=o(n^2)$.

Bai, He, Ma, and Yan~\cite{BHMJ} studied the extremal function
\[
\phi_r(n,k)=\min\{\nu_r(G): |V(G)|=n,\ e(G)=t_{r-1}(n)+k\}.
\]
They proved that $\phi_r(n,k)=(1-o(1))k$ if and only if
$k=o(n^2)$, and determined the common second-order asymptotics
of the integer and fractional packing numbers when
$n=o(k)$ and $k=o(n^2)$. Their work concerns unrestricted graphs,
whereas we investigate how a forbidden-clique condition improves
the optimal packing coefficient.

Forbidding a larger clique can increase the coefficient. Gy\H{o}ri
and Keszegh~\cite{GK} proved that every $K_4$-free graph with
$t_2(n)+k$ edges contains at least $k$ edge-disjoint triangles,
settling a conjecture recorded by Gy\H{o}ri~\cite[p.~360, Conjecture~4]{Gyori1992}.
Liu and Zeng~\cite{LZ} recently extended this result to every fixed
$r\ge3$: for all sufficiently large $n$, every $K_{r+1}$-free graph
with $t_{r-1}(n)+k$ edges contains at least $k$ pairwise edge-disjoint
copies of $K_r$. They also characterize the equality cases.
This is also related to clique decompositions in which each clique
$C$ has cost $|V(C)|-1$. Balogh, He, Krueger, Nguyen, and Wigal proved
the asymptotic integer cost bound $n^2/4+o(n^2)$ for arbitrary
graphs~\cite[Theorem~1.4]{BHKNW}. Their results on fractional edge
decompositions~\cite[Theorems~1.7 and~1.8(ii)]{BHKNW} show that
every graph has a fractional edge decomposition of cost at most
$t_2(n)$. In a $K_4$-free graph, only edges and triangles can occur
in such a decomposition. If $W$ denotes the total weight assigned
to triangles, then the total cost is $e(G)-W$. Thus
$W\ge e(G)-t_2(n)$, and the triangle weights form a fractional
triangle packing. Consequently,
$\nu_3^*(G)\ge e(G)-t_2(n)$.
In a recent preprint, Zhang~\cite{Zhang} proved the exact
integer cost bound $t_2(n)$ for all sufficiently large graphs.
A related line of work concerns fractional clique decompositions,
in which every edge receives total weight exactly one. Longbrake et al.
\cite{LSYZ} establish a Dirac-type sufficient condition for
fractional $K_r$-decompositions of balanced $r$-partite graphs,
while Joos and Smith~\cite{JS} obtain near-optimal probability
thresholds for fractional clique decompositions in random
hypergraphs. These decomposition problems differ from the fractional
packing problem considered here: the total weight on each edge
is required only to be at most one.

However, these results do not determine the optimal dependence of
fixed-size clique-packing bounds on the order of the forbidden clique.

For integers $r\ge3$ and $M\ge r$, define
\[
\gamma_{r,M}=\min\left\{1,\frac{2(M-1)}{r(M-r+1)}\right\}.
\]
The following theorem determines the best coefficient for every
forbidden-clique order. The bound is finite and has no asymptotic error.

\begin{theorem}\label{thm:clique-coefficient}
For all integers $r\ge3$ and $M\ge r$, every $K_{M+1}$-free graph $G$
on $n\ge1$ vertices satisfies
\begin{equation}\label{eq:main-bound}
\nu_r^*(G)\ge\gamma_{r,M}\bigl(e(G)-\rho_{r-1}n^2\bigr).
\end{equation}
For each pair $(r,M)$, the coefficient $\gamma_{r,M}$ is best possible.
\end{theorem}

The coefficient equals $1$ when $M=r$ or $M=r+1$. It is strictly
decreasing for $M\ge r+1$ and tends to $2/r$ as $M\to\infty$.
In terms of a forbidden $K_m$, it is
$\min\{1,2(m-2)/(r(m-r))\}$ for $m\ge r+1$.
The assertion concerns the largest constant valid at all edge
densities. It does not determine the extremal packing number separately
at each density.

For fixed forbidden subgraphs, the same coefficient governs integer
packings asymptotically. Here the error is relative to the excess,
including excesses much smaller than $n^2$.

\begin{corollary}
\label{cor:chromatic-coefficient}
Fix an integer $r\ge3$ and a graph $H$ with $q=\chi(H)\ge r+1$.
Put
\[
\beta_{r,H}=\min\left\{1,\frac{2(q-2)}{r(q-r)}\right\}.
\]
For every $\varepsilon>0$, there exists $n_0$ such that every
$H$-free graph $G$ on $n\ge n_0$ vertices with
$k=e(G)-t_{r-1}(n)>0$ satisfies
$\nu_r(G)\ge(\beta_{r,H}-\varepsilon)k$.
The coefficient $\beta_{r,H}$ is best possible.
\end{corollary}

The chromatic number also characterizes exactly when the
asymptotic coefficient is $1$. This statement includes graphs $H$
with $\chi(H)\le r$, for which positive excesses need not occur
for all sufficiently large orders.

\begin{corollary}\label{thm:chromatic}
Fix $r\ge3$ and a graph $H$. The following property holds if and only
if $\chi(H)\le r+2$. For every $\varepsilon>0$ there is $n_0$ such
that every $H$-free graph $G$ on $n\ge n_0$ vertices with
$k=e(G)-t_{r-1}(n)>0$ satisfies $\nu_r(G)\ge(1-\varepsilon)k$.
\end{corollary}

\begin{remark}
For every fixed $M$, the coefficient in~\eqref{eq:main-bound}
strictly exceeds the unrestricted coefficient in~\eqref{eq:bw}.
When $M\ge r+1$, equality is attained by complete $M$-partite
graphs with equal part sizes; for $M=r$, it is attained by an
unbalanced complete $r$-partite graph. For the special case $H=K_{r+1}$, Liu and Zeng~\cite{LZ}
prove the stronger integer bound
$\nu_r(G)\ge e(G)-t_{r-1}(n)$ for all sufficiently large $n$,
and characterize the equality cases. The fractional theorem does
not, however, eliminate the integer-rounding error or the bounded
difference between $\rho_{r-1}n^2$ and $t_{r-1}(n)$.
\end{remark}

The main difficulty is controlling the number of parts after
symmetrization. For a fixed fractional cover, minimizing its excess
quadratic form over the simplex of part proportions yields, at a
negative minimum, a positive semidefinite matrix with a positive
kernel vector. A Gram representation gives unit vectors whose
pairwise inner products have a common lower bound. We show that
some $r$ of these vectors have a sum satisfying the norm bound
in~\eqref{eq:main-bound}.

From Tinarrage's sharp approximate Carath\'eodory theorem
\cite{Tinarrage}, we derive a sampling lemma with prescribed
inclusion probabilities. The angle condition lets us choose
$r-1$ distinct indices whose vector sum has mean zero.
Adding an index uniformly from the complement gives the required
bound for $r$ indices. In the graph application, this upper bound
on a subset sum contradicts the lower bound imposed by the cover.

Section~\ref{sec:vectors} proves the vector inequality.
Section~\ref{sec:packings} applies it to fractional covers and derives
the integer consequences. Section~\ref{sec:concluding} discusses
the remaining questions.

\section{Subset sums of unit vectors}\label{sec:vectors}

After symmetrization, the problem reduces to a quadratic inequality
in the part sizes. A negative minimum yields a Gram matrix with a
positive kernel vector. We prove a subset-sum bound that rules out
such a minimum.

We use Tinarrage's sharp approximate Carath\'eodory theorem and its
bounded-coordinate version~\cite[Theorem~3 and Lemma~7]{Tinarrage}.
To apply these results to $s$ points in an arbitrary Euclidean space,
we first reduce the dimension to $s-1$. Let $c$ be the orthogonal
projection of the origin onto their affine hull. Then
$\|v_i-c\|^2=\|v_i\|^2-\|c\|^2\le1$, and translation by $c$
preserves every vector combination whose coefficients sum to zero.
We may therefore embed the translated affine hull in
$\mathbb R^{s-1}$. Degenerate configurations follow by continuity:
contract the points slightly toward the origin and perturb them
to be affinely independent while keeping them in the unit ball.
With the convex coefficients fixed, the minimum error over the
finitely many index choices is continuous in the points.

The following lemma allows us to prescribe the mean of the sampled
sum by fixing the inclusion probabilities of the indices.

\begin{lemma}
\label{lem:prescribed-marginals}
Let $v_1,\ldots,v_s$ lie in a Euclidean unit ball, where $s\ge2$.
If $0\le\pi_i\le1$ and $\sum_i\pi_i=k$ is an integer, then there is a
random $k$-element subset $I\subseteq[s]$ such that
\[
\mathbb P(i\in I)=\pi_i,\qquad
\mathbb E\left\|\sum_{i\in I}v_i-\sum_i\pi_i v_i\right\|^2
\le\frac{k(s-k)}{s-1}.
\]
\end{lemma}

\begin{proof}
We maximize the minimum variance over all admissible inclusion
probability vectors. A boundary maximizer is handled by induction
on the number of indices; at an interior maximizer, a supporting
inequality reduces the problem to Tinarrage's lemma.

Write $V=[v_1\ \cdots\ v_s]$, $G=V^{\mathsf T}V$, and
$\Delta=\{x\in[0,1]^s:\sum_i x_i=k\}$. This polytope is the convex
hull of the indicator vectors $\mathbf1_I$ of the $k$-element subsets.
For $x\in\Delta$, let $f(x)$ be the minimum of
$\sum_I p_I\|V\mathbf1_I\|^2$ over all probability distributions
$(p_I)$ with $\sum_I p_I\mathbf1_I=x$. The function $f$ is finite,
continuous, convex, and piecewise affine. Thus $F(x)=f(x)-\|Vx\|^2$ is the minimum variance among
distributions with inclusion probabilities $x$. Indeed, every
feasible distribution satisfies $\mathbb E[V\mathbf1_I]=Vx$, and
therefore
$\mathbb E\|V\mathbf1_I-Vx\|^2
=\mathbb E\|V\mathbf1_I\|^2-\|Vx\|^2$.
The feasible distributions form a nonempty compact polytope,
so the minimum is attained.

We prove the bound by induction on $s$. The cases $k=0$ and $k=s$
are deterministic; the case $s=1$ is also immediate. Thus assume
$1\le k<s$. If $s=2$, then $k=1$, and choosing index $i$ with
probability $x_i$ gives variance
$x_1x_2\|v_1-v_2\|^2\le4x_1x_2\le1$.
Now let $s\ge3$, and choose a maximizer $x^*$ of $F$ on $\Delta$.
If $x_i^*=0$ or $1$, the corresponding indicator is deterministic
and can be removed from the centered sum. The induction hypothesis
then applies with parameters $(s-1,k)$ or $(s-1,k-1)$.
The respective differences between the target bound and the
induction bounds are
\[
\frac{k(k-1)}{(s-1)(s-2)}
\quad\text{and}\quad
\frac{(s-k)(s-k-1)}{(s-1)(s-2)},
\]
both nonnegative.

Suppose next that $x^*$ lies in the relative interior of $\Delta$.
Let $p$ be a subgradient of $f$ on the affine hyperplane
$\sum_i x_i=k$. For every $h$ with $\sum_i h_i=0$ and all
sufficiently small $t$ of either sign, the supporting inequality gives
\[
F(x^*+th)\ge F(x^*)+t(p-2Gx^*)\cdot h-t^2h^{\mathsf T}Gh.
\]
Since $x^*$ maximizes $F$, we have $F(x^*+th)\le F(x^*)$
for both signs of $t$. Thus the linear term vanishes:
$(p-2Gx^*)\cdot h=0$. Consequently, for every vertex
$y=\mathbf1_I$ of $\Delta$,
\[
\|Vy\|^2=f(y)\ge f(x^*)+2(Gx^*)\cdot(y-x^*),
\qquad \|V(y-x^*)\|^2\ge F(x^*).
\]
An optimal distribution defining $f(x^*)$ has average squared
error $F(x^*)$. Hence the minimum squared error over the vertices
is exactly $F(x^*)$. Apply~\cite[Lemma~7]{Tinarrage} with
$d=s-1$ and convex coefficients $x^*/k$. Each coefficient is at
most $1/k$, so the lemma gives a $k$-element set $I$ such that
\[
\left\|\frac1k V\mathbf1_I-\frac1k Vx^*\right\|^2
\le \frac{s-k}{k(s-1)}.
\]
Multiplying by $k^2$ gives a vertex whose squared error is at most
$k(s-k)/(s-1)$. This bounds the maximum of $F$ and proves the lemma.
\end{proof}

Tinarrage's approximate Carath\'eodory theorem allows repeated
indices. The angle condition arising from the covering problem lets
us remove repetitions without increasing the norm.

\begin{lemma}
\label{lem:distinct-indices}
Let $N\ge r\ge3$, and let $v_1,\ldots,v_N$ be unit vectors with
$\langle v_i,v_j\rangle\ge-1/(r-2)$ for $i\ne j$.
Suppose that $\sum_i a_i v_i=0$ for positive weights $a_i$.
For any $r$-term sum with repetitions, there is a sum of $r$ vectors
with distinct indices and no larger norm.
\end{lemma}

\begin{proof}
Let $S$ denote the sum, and suppose that index $i$ occurs
$m\ge2$ times. Then
$\langle S,v_i\rangle\ge m-(r-m)/(r-2)\ge1$, so
$U=S-v_i$ satisfies $\|U\|^2\le\|S\|^2-1$.
If index $j$ occurs $n_j\ge1$ times in $U$, then
$\langle U,v_j\rangle\ge n_j-(r-1-n_j)/(r-2)\ge0$.
The positive relation $\sum_j a_jv_j=0$ therefore guarantees an
unused index $k$ with $\langle U,v_k\rangle\le0$; such an index
exists because $U$ has only $r-1$ terms. Replacing the deleted
occurrence by $v_k$ does not increase the squared norm and increases
the number of distinct indices used. Repeating this step removes
all repetitions and proves the lemma.
\end{proof}

For $N\ge2r$, we have $r\le\lfloor N/2\rfloor$; when
$N=2r-1$, the identity
$r(N-r)=(r-1)(N-r+1)$ gives the same bound.
Applying Tinarrage's theorem~\cite[Theorem~3]{Tinarrage} in
dimension $N-1$, followed by Lemma~\ref{lem:distinct-indices},
yields the bound below for $N\ge2r-1$. This alternative argument is not used later.
Sampling $r-1$ indices first gives a single proof for all
$N\ge r+1$.

\begin{lemma}\label{lem:positive-kernel}
Let $N\ge r+1$ and $r\ge3$. Suppose that $v_1,\ldots,v_N$ are unit
vectors, $\langle v_i,v_j\rangle\ge-1/(r-2)$ for $i\ne j$, and
$\sum_i a_i v_i=0$ for positive weights $a_i$. Then some
$r$-element subset $I\subseteq[N]$ satisfies
\[
  \left\|\sum_{i\in I}v_i\right\|^2\le\frac{r(N-r)}{N-1}.
\]
\end{lemma}

\begin{proof}
Normalize $\sum_i a_i=1$ and put $d=N-r$. Taking the inner product
of $\sum_j a_jv_j=0$ with $v_i$ gives $0\ge a_i-(1-a_i)/(r-2)$,
so $(r-1)a_i\le1$. The numbers $(r-1)a_i$ belong to $[0,1]$ and
sum to $r-1$. Apply Lemma~\ref{lem:prescribed-marginals}
with these inclusion probabilities to obtain a random
$(r-1)$-set $S$. Its sum $U=\sum_{i\in S}v_i$ has mean zero, since
$\mathbb E U=(r-1)\sum_i a_iv_i=0$, and
satisfies $\mathbb E\|U\|^2\le (r-1)(d+1)/(N-1)$.

Choose $J$ uniformly from $[N]\setminus S$ and write $T=\sum_i v_i$.
Since $\mathbb E[v_J\mid S]=(T-U)/(d+1)$ and $\mathbb E U=0$,
\[
  \mathbb E\langle U,v_J\rangle
  =\frac{\mathbb E\langle U,T-U\rangle}{d+1}
  =-\frac{\mathbb E\|U\|^2}{d+1}.
\]
Expanding the squared norm and using $d\ge1$, we obtain
\[
  \mathbb E\|U+v_J\|^2
  =1+\frac{d-1}{d+1}\,\mathbb E\|U\|^2
  \le 1+\frac{(r-1)(d-1)}{N-1}
  =\frac{rd}{N-1}.
\]
Some choice of $I=S\cup\{J\}$ therefore satisfies the bound.
\end{proof}

\section{Clique packings}\label{sec:packings}

We first use fractional symmetrization to reduce the graph problem
to complete multipartite graphs with at most $M$ parts. The vector
inequality then rules out a negative minimum of the covering
quadratic form. Finally, we transfer the fractional bound to
integer packings in graphs with a fixed forbidden subgraph.

A fractional $K_r$-cover assigns nonnegative weights to the edges
so that the weights on every copy of $K_r$ sum to at least one.
By linear programming duality, the minimum total weight of such a
cover equals $\nu_r^*(G)$.
Balogh and Wigal's symmetrization has the following consequence with
a forbidden clique retained.

\begin{corollary}
\label{cor:multipartite-reduction}
Let $r\ge3$ and $s\ge1$ be integers, let $c$ be real, and let $G$ be a $K_{s+1}$-free graph on $n$ vertices. There is a complete multipartite graph $H$ on the same vertex set, with at most $s$ nonempty parts, such that
\[
\nu_r^*(H)+ce(H)\le\nu_r^*(G)+ce(G).
\]
\end{corollary}

\begin{proof}
For the objective $\nu_r^*+ce$, iterate the clone-class
replacements from~\cite[Lemma~2.1]{BW}, as in the proof of
\cite[Lemma~2.2]{BW}. Each step merges two distinct clone classes
whose union is independent and leaves the other classes unchanged,
so the number of classes decreases. The operation preserves
$K_{s+1}$-freeness: any clique meeting the merged class uses at most
one of its vertices and can be mapped to a clique of the same size
in the preceding graph by replacing that vertex with a source vertex.
The process continues until every pair of distinct classes is
completely joined, producing a complete multipartite graph with at
most $s$ nonempty parts.
\end{proof}

We also use the following formula, given in~\cite[proof of Claim~3.2]{BW}.

\begin{lemma}
\label{lem:multipartite-packing}
If $n_1\ge\cdots\ge n_r>0$, then
\[
\nu_r^*(K_{n_1,\ldots,n_r})=n_{r-1}n_r.
\]
A complete multipartite graph with fewer than $r$ nonempty parts has fractional $K_r$-packing number zero.
\end{lemma}

\begin{proof}
Every $K_r$ in a complete $r$-partite graph uses an edge between
the two smallest parts. Summing the packing constraints over the
$n_{r-1}n_r$ edges between these parts gives the upper bound
$n_{r-1}n_r$. For the reverse inequality, assign weight
$1/(n_1\cdots n_{r-2})$ to every $K_r$. An edge between parts
$i$ and $j$ then receives total weight
$n_{r-1}n_r/(n_i n_j)\le1$, while the total packing weight is
$n_{r-1}n_r$. If there are fewer than $r$ parts, the graph contains
no $K_r$, so its fractional packing number is zero.
\end{proof}

The case of at most $r$ parts has a finite linear bound. We record
this consequence of the two known tools to handle small supports.

\begin{proposition}
\label{prop:finite-turan}
For $r\ge3$, every $K_{r+1}$-free graph $G$ on $n\ge1$ vertices satisfies $\nu_r^*(G)\ge e(G)-t_{r-1}(n)$.
\end{proposition}

\begin{proof}
Apply Corollary~\ref{cor:multipartite-reduction} with $c=-1$,
and let $H$ be the resulting complete multipartite graph.
If $H$ has at most $r-1$ parts, then
$e(H)\le t_{r-1}(n)$ and $\nu_r^*(H)=0$.
Otherwise, $H$ has exactly $r$ parts, of sizes
$n_1\ge\cdots\ge n_r$. Merge the two smallest parts.
By Lemma~\ref{lem:multipartite-packing}, the resulting
$(r-1)$-partite graph has
$e(H)-n_{r-1}n_r=e(H)-\nu_r^*(H)$ edges, which is at most
$t_{r-1}(n)$. Thus, in either case,
$\nu_r^*(H)-e(H)\ge-t_{r-1}(n)$.
The reduction transfers this inequality to $G$.
\end{proof}

We apply the vector bound to a minimum of a fixed covering quadratic
form. Working on the closed simplex also takes care of parts whose
weights vanish at the minimum.

\begin{proof}[Proof of Theorem~\ref{thm:clique-coefficient}]
We first prove the lower bound and then establish sharpness.
By symmetrization, it suffices to consider a fixed fractional cover
on a complete multipartite graph. We minimize its excess quadratic
form over the simplex of part proportions and use the vector
inequality to rule out a negative minimum.

\smallskip
\noindent\emph{Reduction and minimization.}
The case $M=r$ follows from Proposition~\ref{prop:finite-turan},
since $t_{r-1}(n)\le\rho_{r-1}n^2$.
Let $M\ge r+1$ and put $C=2(M-1)/(r(M-r+1))\le1$.
By Corollary~\ref{cor:multipartite-reduction}, it suffices to prove
the bound for complete multipartite graphs with at most $M$ parts.
Normalize the part sizes $a_i$ so that they sum to one.
Average a fractional cover over the permutations within each part.
This preserves its total weight and all cover constraints, because these
permutations are graph automorphisms and the feasible region is convex.
The averaged cover is described by nonnegative weights $w_{ij}$,
one for each pair of parts, such that the sum of the weights over
every $r$-element set of parts is at least one. Its normalized cost
is $\sum_{i<j}w_{ij}a_i a_j$, and the normalized excess is
$k=(1/(r-1)-\sum_i a_i^2)/2$.

Suppose, for a contradiction, that a cover violates the required
bound. Fix its weights $w_{ij}$ and minimize
$F(a)=\sum_{i<j}w_{ij}a_i a_j-Ck$ over the closed simplex.
Write $F(a)=a^{\mathsf T}Aa$, where
\[
A_{ii}=\frac{C(r-2)}{2(r-1)},\qquad
A_{ij}=\frac{w_{ij}-C/(r-1)}2\quad(i\ne j).
\]
Let $\eta<0$ be the minimum and restrict to the positive support of
a minimizer $a$, of size $N\le M$. Restrict the weights and the matrix
$A$ to these indices and relabel them by $[N]$. All cover constraints on $r$-element subsets of the support are
inherited from the original cover.
Suppose that $N\le r$. If $k>0$, Proposition~\ref{prop:finite-turan}
implies that the cover cost is at least $k$. To see this first for
rational part proportions, choose a common denominator $m$ and
consider the corresponding blow-up, whose part sizes are $ma_i$.
The restricted weights define a fractional $K_r$-cover of this graph.
Applying Proposition~\ref{prop:finite-turan}, dividing by $m^2$,
and using $t_{r-1}(m)\le\rho_{r-1}m^2$ gives the claimed bound.
Continuity extends it to real part proportions. If $k\le0$,
the nonnegativity of the cover cost suffices. Since $C\le1$,
either case contradicts $F(a)<0$. Hence $N\ge r+1$.

\smallskip
\noindent\emph{From a negative minimum to unit vectors.}
All coordinates of $a$ on the restricted support are positive.
Thus $a+th$ remains in the simplex for every $h$ with
$\sum_i h_i=0$ and all sufficiently small $|t|$.
Expanding $F(a+th)$ and using minimality for both signs of $t$
gives $h^{\mathsf T}Aa=0$ and $h^{\mathsf T}Ah\ge0$.
Hence $Aa$ is a multiple of $\mathbf1$; multiplying by
$a^{\mathsf T}$ and using $\sum_i a_i=1$ yields
$Aa=\eta\mathbf1$.
Set $B=A-\eta J$, where $J$ is the all-ones matrix.
Then $Ba=0$. Writing $z=ta+h$ with $\sum_i h_i=0$, we obtain
$z^{\mathsf T}Bz=h^{\mathsf T}Ah\ge0$.
Thus $B$ is positive semidefinite.
Its diagonal entries equal
$\beta=C(r-2)/(2(r-1))-\eta>0$.
Represent $B/\beta$ as the Gram matrix of unit vectors $v_i$.
Since $Ba=0$, the squared norm of $\sum_i a_iv_i$ is
$a^{\mathsf T}Ba/\beta=0$. Thus $\sum_i a_iv_i=0$, and
\[
B_{ij}+\frac{\beta}{r-2}
 =\frac{w_{ij}}2-\left(1+\frac1{r-2}\right)\eta\ge0.
\]
Their pairwise inner products are therefore at least $-1/(r-2)$.
\smallskip
\noindent\emph{The contradiction.}
For every $r$-element set $I$, expanding the Gram sum gives
\[
\begin{aligned}
\beta\left\|\sum_{i\in I}v_i\right\|^2
&=\sum_{i,j\in I}B_{ij}\\
&=\sum_{\substack{i<j\\i,j\in I}}w_{ij}
 -\frac{Cr}{2(r-1)}-r^2\eta.
\end{aligned}
\]
The cover constraint on $I$ therefore yields
\begin{equation}
\beta\left\|\sum_{i\in I}v_i\right\|^2
 \ge1-\frac{Cr}{2(r-1)}-r^2\eta.
\label{eq:cover-gram-constraint}
\end{equation}
By Lemma~\ref{lem:positive-kernel}, some $r$-element set has
squared norm at most $r(N-r)/(N-1)\le L$, where
$L=r(M-r)/(M-1)$. The last inequality follows from $N\le M$
and the fact that $(N-r)/(N-1)=1-(r-1)/(N-1)$ is increasing in $N$.
Substituting this upper bound into
\eqref{eq:cover-gram-constraint} gives
\[
1\le \beta L+\frac{Cr}{2(r-1)}+r^2\eta
=\frac{C(r+(r-2)L)}{2(r-1)}+(r^2-L)\eta.
\]
By the definitions of $C$ and $L$,
$C(r+(r-2)L)/(2(r-1))=1$. Since $L<r^2$ and $\eta<0$, we obtain
$1\le1+(r^2-L)\eta<1$, a contradiction.
Linear programming duality and the multipartite reduction now
give the desired bound.

\smallskip
\noindent\emph{Sharpness.}
For $M\ge r+1$, let $G$ be the complete $M$-partite graph
with every part of size $t$. Each edge lies in
$\binom{M-2}{r-2}t^{r-2}$ copies of $K_r$, so assigning the
reciprocal of this number to each copy saturates every edge.
This gives the optimal packing weight $e(G)/\binom r2$.
The fractional packing number and the excess are, respectively,
\[
\frac{M(M-1)t^2}{r(r-1)}
\quad\text{and}\quad
\frac{M(M-r+1)t^2}{2(r-1)}.
\]
Their ratio is $\gamma_{r,M}$.
For $M=r$, take a complete $r$-partite graph with two parts
of size $t$ and the other $r-2$ parts of size $2t$.
Merging the two smaller parts removes $t^2$ edges and gives
the balanced complete $(r-1)$-partite graph on $2(r-1)t$
vertices, whose edge count is $\rho_{r-1}(2(r-1)t)^2$.
Thus the excess is $t^2$. Every $K_r$ uses one of the $t^2$
edges between the two smaller parts, so the fractional packing number is at most $t^2$. Assigning weight $(2t)^{-(r-2)}$ to each $K_r$
attains this bound. Hence the coefficient is sharp for every
pair $(r,M)$.
\end{proof}

To obtain an error relative to the excess, we combine the
dense-range transfer with the stronger estimate available
for small excesses.

\begin{proof}[Proof of Corollary~\ref{cor:chromatic-coefficient}]
By the theorem of Erd\H{o}s, Frankl, and
R\"odl~\cite[Theorem~1.5]{EFR}, every $n$-vertex $H$-free graph
can be made $K_q$-free by deleting at most $a(n)=o(n^2)$ edges.
Let $G'$ be the resulting graph. Let $b(n)=o(n^2)$ be the
uniform error in the Haxell--R\"odl theorem~\cite{HR}, as stated
in Theorem~1.7 of Balogh and Wigal~\cite{BW}, and put
$g(n)=\rho_{r-1}n^2-t_{r-1}(n)=O_r(1)$.
Writing $\beta=\beta_{r,H}$, Theorem~\ref{thm:clique-coefficient}
with $M=q-1$ gives
\[
\nu_r(G)\ge\nu_r(G')\ge\nu_r^*(G')-b(n)
\ge\beta k-\beta\bigl(a(n)+g(n)\bigr)-b(n).
\]
These estimates hold uniformly over $G$. Choose
$h(n)=n\sqrt{a(n)+b(n)+g(n)+1}$.
Then $h(n)=o(n^2)$ and
$a(n)+b(n)+g(n)+1=o(h(n))$.
For $k\ge h(n)$, we obtain
$\nu_r(G)\ge(\beta-o(1))k$.
For $0<k<h(n)$, Gy\H{o}ri's small-excess theorem, stated for
all fixed $r\ge3$ in Theorem~1.6 of Balogh and Wigal~\cite{BW}, gives
$\nu_r(G)\ge k-O_r(k^2/n^2)
\ge(1-O_r(h(n)/n^2))k$.
Since $\beta\le1$, both ranges yield the required bound
uniformly over $G$ and $k$.

For sharpness when $q\ge r+2$, take
$G=T_{q-1}(n)$, with $n$ divisible by $q-1$.
These graphs are $H$-free. The elementary bound
$\nu_r(G)\le e(G)/\binom r2$ gives
\[
\limsup_{n\to\infty}
\frac{\nu_r(T_{q-1}(n))}{e(T_{q-1}(n))-t_{r-1}(n)}
\le\frac{2(q-2)}{r(q-r)}=\beta.
\]
When $q=r+1$, take the complete $r$-partite graph with $r-2$
parts of size $2t$ and two parts of size $t$. This graph is
$H$-free. Merging the two smaller parts gives
$T_{r-1}(2(r-1)t)$, so the excess is $t^2$. Every $K_r$ uses
an edge between the two smaller parts, and there are only $t^2$
such edges. Hence its packing number is at most $t^2$, showing
that the coefficient cannot exceed $1$.
\end{proof}

\begin{proof}[Proof of Corollary~\ref{thm:chromatic}]
Suppose first that $\chi(H)\le r+2$. Choose a complete
$(r+2)$-partite graph $F$ containing $H$, with every part nonempty.
Every $H$-free graph is then $F$-free. Since $\chi(F)=r+2$,
Corollary~\ref{cor:chromatic-coefficient} applied to $F$ gives
the coefficient-$1$ property. Conversely, if
$q=\chi(H)\ge r+3$, that corollary gives the best coefficient
$2(q-2)/(r(q-r))<1$, so the coefficient-$1$ property fails.
This proves the equivalence, including the case $\chi(H)\le r$,
which is not covered by the coefficient corollary.
\end{proof}

\section{Concluding remarks}\label{sec:concluding}

Theorem~\ref{thm:clique-coefficient} determines the optimal uniform
linear coefficient for every pair of clique orders. The corresponding
fixed-density problem is more refined. For $K_{r+1}$-free graphs,
the optimal uniform coefficient is $1$, although the optimal density
bound exceeds the excess on part of the density range. A separate paper~\cite{HeDensities}
determines that density function and a finite integer form for
$K_4$-free triangle packings. It also treats the complete triangle
density curve in $K_5$-free graphs. These results are not used in the
proofs here. An earlier publicly circulated version of that paper
included the cases $M=r+1$ and $r=3$ of our coefficient theorem,
together with the coefficient-$1$ chromatic criterion. Those results
are incorporated into the general treatment here; the revised
companion paper concentrates on density bounds and finite refinements.
For larger forbidden cliques, determining the entire density function
remains open.

The normalization by $\rho_{r-1}n^2$ cannot always be replaced by
$t_{r-1}(n)$ in the finite theorem. The graph $K_5$, which is
$K_6$-free, has $\nu_4^*(K_5)=5/3$, whereas
$e(K_5)-t_3(5)=2$. For triangles, the same graph has
$\nu_3^*(K_5)=10/3<(8/9)(e(K_5)-t_2(5))$.
Thus even when the uniform coefficient is known, finite corrections
require additional information.

The equality constructions prove sharpness but do not classify
all equality cases. A stability theorem would require a description
of both near-extremal fractional covers and part sizes.
Corollary~\ref{cor:chromatic-coefficient} leaves an $o(k)$ error
in the integer problem; determining the optimal error for individual
forbidden graphs remains open. The vector inequality also suggests
studying which angle and dependence conditions yield sharp
subset-sum bounds in other fractional packing problems.

\paragraph{Declaration of AI assistance.}
OpenAI Codex was used to assist with mathematical discussions, proof review, and language editing during the preparation of this paper. The authors take full responsibility for the content and correctness of the paper.

\end{document}